\documentclass[11pt]{amsart}

\usepackage{amssymb, amsmath, epsfig, graphics, amscd, ifthen}
\usepackage[all, cmtip]{xy}

\newcommand{\pics}{true}

\newtheorem{theorem}{Theorem}[subsection]

\newtheorem{proposition}[theorem]{Proposition}
\newtheorem{corollary}[theorem]{Corollary}
\newtheorem{lemma}[theorem]{Lemma}
\newenvironment{pf}{\paragraph{\sc Proof}}{\hspace*{\fill}$\Box$\par\medskip}

\newtheorem{example}[theorem]{Example}

\newtheorem{remark}[theorem]{Remark}

\newcommand{\aut}{{\mathop{\rm Aut}\nolimits}}
\newcommand{\spec}{\mathop{\rm Spec}\nolimits}

\renewcommand{\mod}{\mathop{\rm mod}\nolimits}
\newcommand{\Amod}{A{\rm -mod}}
\newcommand{\nilp}{{\rm nilp}}
\newcommand{\Anilp}{A{\rm -nilp}}
\newcommand{\AMod}{A{\rm -Mod}}
\newcommand{\rep}{{\rm rep}}

\newcommand{\pic}{\mathop{\rm Pic}\nolimits}

\newcommand{\Coh}{\mathop{\rm Coh}\nolimits}

\renewcommand{\phi}{\varphi}

\newcommand{\Hom}{\mathop{\rm Hom}\nolimits}

\newcommand{\Ext}{\mathop{\rm Ext}\nolimits}
\newcommand{\End}{\mathop{\rm End}\nolimits}

\newcommand{\im}{\mathop{\rm Im}\nolimits}

\newcommand{\Stab}{\mathop{\rm Stab}\nolimits}
\newcommand{\DAut}{\mathop{\rm DAut}\nolimits}

\newcommand{\SL}{{\rm SL}}
\newcommand{\GL}{{\rm GL}}

\newcommand{\Hilb}{{\rm Hilb}}

\newcommand{\Tr}{{\rm Tr}}

\renewcommand{\Box}{\square}
\renewcommand{\tilde}{\widetilde}

\newcommand{\F}{\mathcal F}

\newcommand{\M}{\mathcal M}

\newcommand{\D}{\mathcal D}
\newcommand{\N}{\mathbb N}

\newcommand{\I}{\mathcal I}

\newcommand{\sA}{\mathcal A}
\newcommand{\sP}{\mathcal P}

\newcommand{\E}{\mathcal E}
\newcommand{\R}{\mathbb R}

\newcommand{\V}{\mathcal V}

\newcommand{\C}{\mathbb C}

\newcommand{\Z}{\mathbb Z}

\renewcommand{\O}{\mathcal O}
\renewcommand{\L}{\mathcal L}
\newcommand{\PP}{\mathbb P}

\newcommand{\wt}{{\rm wt}}

\newcommand{\q}{{\bf q}}
\renewcommand{\v}{{\bf v}}

\newcommand{\color}[6]{} 

\begin{document}

\begin{center}
{\LARGE Non-commutative Donaldson--Thomas theory

\vspace{0.1in}

and the conifold}

\vspace{0.2in}

{\large Bal\'azs Szendr\H oi}

\vspace{0.1in}

{Mathematical Institute, University of Oxford, UK}

\vspace{0.1in}

{\tt szendroi@maths.ox.ac.uk}

\end{center}

\vspace{0.15in}

\thispagestyle{empty}

{\small
\begin{center} {\sc abstract} \end{center}
{\leftskip=50pt \rightskip=50pt
\noindent Given a quiver algebra~$A$ with relations defined by a superpotential, 
this paper defines a set of invariants of~$A$ counting framed cyclic $A$-modules, 
analogous to rank-1 Donaldson--Thomas invariants of Calabi--Yau threefolds. 
For the special case when~$A$ is the non-commutative crepant resolution of the 
threefold ordinary double point, it is proved using torus localization that 
the invariants count certain pyramid-shaped partition-like 
configurations, or equivalently infinite dimer configurations in 
the square dimer model with a fixed boundary condition. The 
resulting partition function admits an infinite product expansion,
which factorizes into the rank-1 Donaldson--Thomas partition functions
of the commutative crepant resolution of the singularity and its flop. 
The different partition functions are speculatively interpreted as counting stable 
objects in the derived category of $A$-modules under different stability 
conditions; their relationship should then be an instance of wall crossing 
in the space of stability conditions on this triangulated category.
\par}}

\section*{Introduction}

This paper is concerned with phenomena in the enumerative geometry of 
local Calabi--Yau threefolds, mostly via a study of 
the example $X=\O_{\PP^1}(-1, -1)$, a toric Calabi--Yau variety 
often referred to as the resolved conifold. It is known 
that topological string theory associates a function $Z_X$ of two variables~$(q,t)$ 
to~$X$, its topological string partition function, admitting
the following infinite product form: 
\[Z_X(q,t) = M(-q)^2\prod_{k\geq 1} \left(1-(-q)^k e^{-t}\right)^k,
\]
where the MacMahon function $M$ is itself an infinite product
\[M(q)=\prod_{k\geq 1}(1-q^k)^{-k}.
\]
According to~\cite{mnop1}, the function $Z_X$ admits two different Laurent 
series expansions, with coefficients which are enumerative invariants 
associated to~$X$. The expansion in the variables $(q, e^{-t})$ gives
coefficients which are Donaldson--Thomas invariants, defined as the virtual 
number of singular $U(1)$-connections (more precisely
rank-1 torsion free sheaves) on~$X$, evaluated by enumetating certain 
combinatorial arrangements related to 3-dimensional partitions
(compare~\textsection\ref{sec_cdt}). On the other hand, dividing by the MacMahon 
factor, we can expand in variables $(u, e^{-t})$, with 
the new variable~$u$ related to~$q$ via $q=-e^{iu}$. The coefficients in 
this case are the (disconnected) Gromov--Witten invariants of~$X$, counting stable 
maps of (not necessarily connected) curves to~$X$ of arbitrary genus and degree. 

In the language of string theory, 
the variable~$u$ is the string coupling constant. The two expansions are 
in the parameter region (large~$t$, small~$u$), where perturbative string theory
gives a good description and hence the coefficients are 
Gromov--Witten invariants, and 
(large~$t$, large $|u|$), the region where the physical description is
(conjecturally) 
a version of $U(1)$ gauge theory~\cite{inov}, leading to rank-1 Donaldson--Thomas 
invariants. This leads to an obvious question: what happens in parameter regions
where~$t$ is small? 

The parameter~$t$ is geometric: it is the K\"ahler class, measuring the volume 
of the zero-section $\PP^1\cong C\subset X$. Hence the limit $t\to 0$ corresponds
to contracting the zero-section, resulting in the conifold singularity
$Z=\{x_1x_2-x_3x_4=0\}\subset \C^4$. This singular variety is
not known to have sensible enumerative invariants. I consider instead its
non-commutative resolution~\cite{vdb2}, a homologically smooth 
non-commutative Calabi--Yau algebra~$A$, and show that this algebra gives
rise to a set of enumerative invariants analogous to the Donaldson--Thomas 
invariants of~$X$. These invariants are combinatorial, like those
of~$X$, this time given by counting
finite subsets of a certain rectangular arrangement of two-coloured stones, 
or equivalently, infinite dimer configurations on the square lattice
with fixed asymptotics. Using a combinatorial result due to Young~\cite{young},
a non-perturbative change of variables leads to 
a partition function with a product expansion
\[Z_A(q,t) = M(-q)^2\prod_{k\geq 1} \left(1-(-q)^k e^{-t}\right)^k\left(1-(-q)^k e^{t}\right)^k, 
\]
closely related to that of $Z_X$, which however 
only admits a Laurent expansion near $e^{-t}=1$ or $t=0$. Thus, 
in this parameter region, a new phenomenon emerges: the function 
$Z_X$ is roughly ``the positive half'' of the full function $Z_A$, which is 
given naturally by non-commutative geometry and dimer models. 

Dimers have been in fashion lately in the high energy 
physics literature, already making
an occurence in~\cite{orv}, and independently in work of Hanany and 
coauthors~\cite{hk, fh} studying gauge/string duality for toric 
Calabi--Yau singularities. However, their direct appearance in the mathematical 
literature, connecting them to enumerative invariants of non-commutative 
Calabi--Yau geometries, seems to be new. It would be interesting to relate
this idea to the other direct derivation of dimer models from toric 
geometry~\cite{fengetal}, based on mirror symmetry and special Lagrangian geometry.

Section~1 studies quiver algebras defined by a superpotential, of which the conifold
algebra~$A$ is an example, proving in~\textsection\ref{sec_obs} 
the basic fact that moduli 
spaces of framed cyclic modules over such algebras (the exact analogues 
of rank-1 torsion-free sheaves, or rather the corresponding quotient sheaves) 
admit a perfect obstruction theory. The obstruction theory is constructed
via natural embeddings of these highly singular moduli spaces into smooth
varieties; the existence of these embeddings may be of independent interest.

Section~2 is about the non-commutative conifold, reducing
the computation of its enumerative invariants to combinatorics and dimer
configurations in~\textsection\ref{sec_pyr}-\textsection\ref{sec_dim}, 
and discussing the product form of the partition function 
in~\textsection\ref{sec_prod}. 

Section~\ref{chap_interp} discusses possible interpretations and generalizations.
The natural space paramet\-rized by the variable $z=e^{-t}$ 
is introduced in~\textsection\ref{sec_br}; this is the space of Bridgeland
stability conditions of a variant of the derived category of coherent sheaves 
on~$X$. The statement that the limit $t\to 0$ 
should be thought of as moving to the non-commutative resolution~$A$
becomes more precise on this space. 
A speculative interpretation of the partition functions $Z_X$, $Z_A$ as counts 
of stable objects in the derived category is also given there. This points towards
an interpretation of the change from $Z_X$ to $Z_A$ as a result
of performing a countably infinite number of wall crossings
in the space of stability conditions.

A curious expansion of the partition function $Z_A$ near small~$t$, involving 
Eisenstein-like sums, is discussed in~\textsection\ref{sec_mod}.
In~\textsection\ref{sec_higher} an obvious combinatorial generalization of the 
partition function is studied, though the enumerative interpretation of 
this generalization is currently unclear. \textsection\ref{sec_glob} presents 
a way of extending the problem studied here to a global situation.  
The concluding~\textsection\ref{sec_orb}
briefly discusses work of Bryan--Young~\cite{by} on some orbifold examples, 
where the partition function exhibits the same ``doubling'' phenomenon when 
moving from the commutative resolution to the non-commutative one.

\section{Non-commutative Donaldson--Thomas theory}

\subsection{Quiver algebra and the superpotential}

Let $Q=\{V, E\}$ be a quiver with vertex set~$V$, edge set~$E$ and head and tail 
maps $h, t\colon E\to V$ prescribing the head and tail of each oriented edge 
(which may in general coincide). Let $\C Q$ be the path algebra of~$Q$, generated
over $\C$ by oriented paths, with multiplication defined by joining
paths. The vector space $\C Q/[\C Q,\ CQ]$ has as basis the set of cycles up to 
cyclic permutation of the arrows.
Given a {\em superpotential}, an element $W\in\C Q/[\C Q, \C Q]$, 
define a two-sided ideal $I_W$ of $\C Q$ by
\[ I_W = \langle\langle \partial_e W \ | \ e\in E \rangle\rangle.
\]
Here represent $W$ as an element in the path algebra consisting of 
a sum of cycles, well-defined up to cyclic permutation. Define formal 
differentiation $\partial_e$ of a cyclic monomial with respect to an 
edge $e\in E$ to be zero if~$e$ does not appear in the monomial;
otherwise, cyclically permute the monomial until the edge~$e$ is in 
the first position, and delete it. For details, 
see e.g.~\cite{bocklandt, ginzburg}.
Let $A=\C Q/I_W$ be the quotient of the path algebra $\C Q$ 
by the ideal $I_W$ of relations defined by the superpotential~$W$.

The algebra $A$ contains a set $\{f_i:i\in V\}$ of orthogonal idempotents based at 
the vertices of $Q$, represented by paths of length~$0$. 
Thus $A$ decomposes as $A=\oplus_{i\in V} Af_i$ into a set $P_i=Af_i$ of projective
left $A$-modules. Note that here and everywhere below, juxtaposition denotes
multiplication in the algebra~$A$, or in an $A$-module, as appropriate.
The algebra~$A$ is generated 
over the commutative idempotent ring $\C[f_i : i\in V]$ by non-commuting variables 
$\{x_e: e\in E\}$ attached to the arrows of the quiver, which satisfy the relations 
obtained by formal differentiation of $W$ with respect to the variables $x_e$. 
As $A$ is a path algebra of a quiver with relations, the category 
$\Amod$ of finite-dimensional left $A$-modules
is equivalent to the category $\rep_{I_W}(Q)$ of finite-dimensional representations 
$(M_i, \phi_e)$ of the quiver satisfying the relations~$I_W$. I will switch 
between the two languages freely without comment. Objects in $\rep_{I_W}(Q)$ have a 
dimension vector $\v=\{v_i=\dim M_i\}\in\Z^V$.

Let $M$ be a cyclic left $A$-module, generated by an element $m\in M$; thus 
$M=Am$. Then $M=\oplus_{i\in V} Am_i$ where $m_i=f_im$. 
A cyclic module generated by a vector $m\in f_kM$ will be called a cyclic module
based at the vertex $k\in V$.

\begin{example}\rm\label{ex_threevar} Consider the quiver $Q=\{V,E\}$ with one
vertex $V=\{*\}$, three loop edges $E=\{x_1,x_2,x_3\}$ with heads and tails at the 
given vertex $*$, and cubic superpotential 
\[W=x_1x_2x_3-x_1x_3x_2.\]
The quiver algebra
\[\C Q = \C\langle x_1, x_2, x_3\rangle\]
is a free non-commutative $\C$-algebra on three generators. 
The ideal $I_W$ in $\C Q$ is generated by the relations
\begin{eqnarray*}
 \partial_{x_1} W & = & x_2x_3-x_3x_2,\\ 
 \partial_{x_2} W & = & x_3x_1-x_1x_3,\\ 
 \partial_{x_3} W & = & x_1x_2-x_2x_1.\\ 
\end{eqnarray*}
Thus
\[ A = \C\langle x_1, x_2, x_3\rangle / \langle x_ix_j-x_jx_i \rangle\cong \C[x_1,x_2,x_3]
\] 
is commutative, the coordinate ring of affine space~$\C^3$. 
The dimension lattice of~$A$ is~$\Z$.
 
A more substantial example, relevant to the rest of the paper, is discussed
in~\textsection\ref{sec_conifoldA}.
\end{example}

\subsection{The moduli space of cyclic modules}
\label{sec_modconst}

Let~$A$ be an algebra defined by a quiver $Q=\{V,E\}$ with 
superpotential~$W$. Fix a vertex $k\in V$, as well as 
vector spaces $\{U_i:i\in V\}$ with dimension 
vector $\v\in \Z^V$. Consider the vector space
\[ S= \left(\prod_{e\in E} \Hom(U_{t(e)}, U_{h(e)})\right) \times U_k=\left\{(\phi_e)_{e\in E}, m \right\}.\] 
Let $S^0\subset S$ be the open subvariety
defined by the condition that the vector $m\in U_k$ generates the
$\C Q$-module $(U_i, \phi_e)$. Let
\[ X = \left\{ \partial_e W =0 \,\, \Big| \,\, e\in E \right\}\subset S^0\]
be the closed subscheme of $S^0$ cut out by the superpotential equations. 

The group $G=\prod_{i\in V} \GL(U_i)$ acts on $S$ by 
\[ (g_i) \circ \left((\phi_e), m\right) = \left((g^{}_{h(e)}\phi^{}_e g_{t(e)}^{-1}), g_km\right),
\] 
where $h, t:E\to V$ are the head and tail maps.
Lift the action to $S\times\C$ by 
\begin{equation}\label{lin} (g_i) \circ ((\phi_e), m, z) = \left((g^{}_{h(e)}\phi^{}_e g_{t(e)}^{-1}), g_km, \chi^{-1}(g_i)z\right),
\end{equation}
where the character $\chi$ of~$G$ is defined by 
$\chi(g_i)=\prod_{i\in V} \det(g_i)$. The following is now standard~\cite{king}: 

\begin{lemma}\label{lem_stab}\begin{enumerate}
\item The action of~$G$ on~$S^0$ is free.
\item For a point $((\phi_e), m, z)\in S\times \C$ with $z\neq 0$, the closure of
the orbit $\overline{G((\phi_e), m, z)}\subset S\times \C$ is disjoint from the 
zero-section $S\times\{0\}$ if and only if $((\phi_e), m)\in S^0$.
\end{enumerate}
\end{lemma}
\begin{pf}
For $(1)$, if $g\in G$ fixes $(\{\phi_e\}, m)\in S^0$, then 
$\ker(g-{\rm id})\subset\oplus_{i\in V} U_i$ is a 
subspace containing $m$ and stable under $\langle \phi_j\rangle$. 
Thus, by definition, $g={\rm id}$. Hence the action of~$G$ on $S^0$ is free. 

For (2), suppose that $((\phi_e), m)$ is not cyclic. Then there are 
decompositions $U_i = V_i\oplus W_i$ with 
\[\phi_e = \left(\begin{array}{cc}* & * \\ 0 & *\end{array}\right) \ \mbox{ and }\ 
m = \left(\begin{array}{c}* \\ 0 \end{array}\right)\]
with respect to this decomposition. For $s\in\C^*$, let
\[g_i(s) = \left(\begin{array}{cc}1 & 0 \\ 0 & {\rm diag}(s^{-1})\end{array}\right),\]
then 
\[ (g_e(s))\circ \left(\left(\begin{array}{cc}* & * \\ 0 & *\end{array}\right), m, z\right) = \left(\left(\begin{array}{cc}* & s* \\ 0 & *\end{array}\right), m, s^Nz\right)
\]
with $N>0$, having a limit at $s=0$ on the zero-section.

Conversely, assume that $((\phi_e), m)$ is cyclic, but
the closure of the orbit $G((\phi_e), m, z)$ intersects the zero-section.
Then there is a one-parameter subgroup $\lambda\colon\C^*\to G$ with
\[\lim_{s\to 0}\lambda(s)\circ((\phi_e), m, z)\in S\times\{0\}.\]
Decompose $U_i=\oplus_{n\in\Z} U_{i,n}$ under the
action of $\C^*$. Since the limit at $s=0$ exists, $m\in \oplus_{n\geq 0} U_{i,n}$, 
and $\phi_e\colon U_{t(e), m}\to\oplus_{n\geq 0} U_{h(e),m+n}$. Thus, since
$((\phi_e), m)$ is cyclic, $U_i=\oplus_{n\geq 0} U_{i,n}$. Hence 
$\chi(\lambda(s))=s^N$, with $N\geq 0$. $N=0$ is impossible, since then the limit
would be contained in the orbit. If $N>0$, then 
\[\lambda(s)\circ ((\phi_e), m, z) = ((\lambda^{}_{h(e)}(s)\phi^{}_e\lambda^{-1}_{t(e)}(s)),\lambda(s)m, s^{-N}z)
\]
does not converge to a point on the zero-section as $s\to 0$, a contradiction. 
\end{pf}

\begin{proposition} \begin{enumerate}\item 
There exists a smooth and quasi-projective 
geometric quotient~$N$ of~$S^0$ by~$G$, containing a closed 
subscheme $\M_{k,\v}\subset N$ which is a quotient of~$X$ by~$G$.
\item The space $\M_{k,\v}$ carries a tautological family $(M_{k,\v}, m_{k,\v})$
of framed cyclic $A$-modules, generated at the vertex~$k$. 
\item The space $\M_{k,\v}$ a fine moduli space; the triple 
$(\M_{k,\v}, M_{k,\v}, m_{k,\v})$ represents the functor of flat 
families of framed $A$-modules (locally free sheaves with $A$-structure
and $A$-generator) over schemes, generated by a section based 
at the vertex~$k$.
\end{enumerate} 
\label{prop_modspace}
\end{proposition}
\begin{pf} To prove (1), following~\cite{king}, regard the affine space $S$ as a 
quasi-projective variety, and consider the linearization~(\ref{lin}) of 
the action of~$G$ on the trivial line bundle $S\times\C$. 
Lemma~\ref{lem_stab}~(2) implies that the semistable locus $S^{\rm ss}$ is exactly 
the subset $S^0$. This also agrees with the stable locus, since all stabilizers
of points in $S_0$ are trivial by Lemma~\ref{lem_stab}~(1).
Thus, by Geometric Invariant Theory, a quasi-projective geometric 
quotient $N=S/\!/_{\chi} G = S^0/G$ exists, and it is smooth because 
stabilizers on $S^0$ are trivial. Since~$X$ is $G$-invariant in $S^0$, 
its quotient by~$G$ exists as a closed subscheme $\M_{k,\v}\subset N$.

For (2), it is enough to observe that, tautologically, $S^0$ carries a 
tautological family of $\C Q$-modules, generated by a section. Over 
$X\subset S^0$, the relations are also satisfied, and hence the family
becomes that of $A$-modules. Being $G$-equivariant, the family and its generator
descend to the space $\M_{k,\v}$. 

Finally, to see (3), let $Z$ be a scheme with a flat family $M_Z$ of $A$-modules,
locally free sheaves with $A$-structure, generated by a section $m\in H^0(Z, M_Z)$.
Taking trivializations of $M_Z$ on an open cover $\{Z_j\}$ of $Z$ gives 
tautological maps $Z_j\to X$, and composing with the projection $X\to\M_{k,\v}$, 
these maps glue to a map $Z\to\M_{k,\v}$ under which the family on $Z$ is a 
pullback by construction. 
This shows that $(\M_{k,\v}, M_{k,\v}, m_{k,\v})$ is indeed a universal family. 
\end{pf}

At a point $[M, m]\in \M_{k,\v}$, let $\bar m \colon A \to M$ be the canonical 
$A$-module surjection given by $a\mapsto am$. Let $I=\ker\bar m$. Since $m\in f_kM$,
where $f_k$ is the idempotent based at the $k$-th vertex, the ideal $I$ decomposes
as a direct sum 
$I=I_k\oplus \bigoplus_{j\neq k} P_j$, into the left $A$-ideal $I_k=If_k$
and the remaining projective modules $P_j = Af_j$ for $j\neq k$. 

\begin{corollary} The Zariski tangent space to $\M_{k,\v}$ at its point 
$[M, m]$ can be identified with the $\C$-vector space $\Hom_A(I_k, M)$. 
\label{cor_tangent}
\end{corollary}
\begin{pf} Each deformation $(M_\eta, m_\eta)$ of $(M, m)$ in $\M_{k,\v}$
is still based at the vertex~$k$. Hence the projectives $P_j$ for $j\neq k$ must 
remain in the kernel of $\bar m_\eta$. With this modification, the 
statement is a standard corollary of Proposition~\ref{prop_modspace} (2);
see e.g.~\cite[Proposition E2.4(ii)]{az}.
\end{pf}

\subsection{Perfect obstruction theory on the moduli space}
\label{sec_obs}

Recall~\cite{bf0} that given a scheme $Z$, a perfect obstruction theory 
$(E, \psi)$ for~$Z$ is morphism
$\psi: E \to L_Z$ in the derived category $\D(Z)$ of quasicoherent $\O_Z$-modules, 
where $L_Z$ is the cotangent complex of $Z$, $E\in\D(Z)$ is a perfect complex of 
amplitude $[-1, 0]$, and $\psi$ induces an isomorphism on $H^0$ and a surjection 
on $H^{-1}$. A symmetric perfect obstruction theory~\cite{bf} is 
a triple $(E, \psi, \theta)$ including also a non-degenerate 
symmetric bilinear form $\theta: E\to E^\vee[1]$.

\begin{theorem} The moduli space $\M_{k,\v}$ carries 
a symmetric perfect obstruction theory.
\label{th_obs}
\end{theorem}
\begin{pf} Recall the superpotential $W\in \C Q/[\C Q, \C Q]$ used to define 
the algebra $A$. Consider the regular function $w=\Tr(W)$ on 
the smooth variety~$S^0$, 
obtained by taking the sum of the traces of the cycles making up~$W$. 
As spelled out for example in~\cite[Proposition 3.8]{segal}, 
\[X=Z(dw)\subset S_0\] is exactly the scheme-theoretic vanishing locus of the 
one-form $dw\in\Omega^1_{S^0}$. Note that~\cite{segal} deals with the case where
there is no cyclic generator, but the proof carries over verbatim: the 
equations defining $X\subset S_0$ are the relations between various linear maps
prescribed by the superpotential~$W$; they do not involve the generator.
The cyclic generator is only used to define the open subset $S_0\subset S$. 

The function $w\in H^0(\O_{S^0})$ is invariant under the $G$-action above, 
and hence descends to a regular function $\bar w\in H^0(\O_N)$ on the smooth 
quotient~$N$. By
naturality, $\M_{k,\v}=Z(d\bar w)$ is the vanishing locus on $N$ of the exact 
one-form~$d\bar w$. By~\cite[Remark 3.12]{behrend}, this defines 
a symmetric perfect obstruction theory on~$\M_{k,\v}$. 
\end{pf}

\begin{remark}\rm Recall that an associative $\C$-algebra $A$ is 
{\em 3-Calabi--Yau}, if for all $M, N\in\AMod$, the category of 
finitely generated $A$-modules, with at least one of $M$, $N$ finite-dimensional,
there exist perfect bifunctorial pairings
\[ \Ext^i_A(M,N) \times \Ext^{3-i}_A(N,M) \to \C\]
between finite-dimensional $\C$-vector spaces. 
For certain choices of superpotentials, it is known that the algebra~$A$ considered 
above is 3-Calabi--Yau~\cite{bocklandt, ginzburg}. 
The Calabi--Yau duality between $\Ext^1$ and $\Ext^2$ is 
morally responsible for the existence of the symmetric 
perfect obstruction theory on the module space of $A$-modules, though as the proof 
shows, it is not a necessary requirement for cyclic modules.

For the case of sheaves on a Calabi--Yau threefold, the construction of the 
symmetric obstruction theory is due to Thomas. The proof in~\cite{th_th} 
uses the equations of the moduli space inside a Grassmannian; this inspired
the proof of Theorem~\ref{th_obs} above. The proof of~\cite{th} uses duality 
explicitly. 
\end{remark}

\subsection{Cyclic modules and finite-codimension ideals}
\label{sec_ideals}

The purpose of this section is to show that, at least set-theoretically and
under the Calabi--Yau assumption, parametrizing cyclic finite-dimensional
$A$-modules up to isomorphism, respecting the cyclic generator, is equivalent to 
parametrizing finite-codimension ideals of~$A$ up to $A$-module isomorphism. 
This is analogous to the relationship between the moduli space of ideal sheaves and
the Hilbert scheme on a projective variety. I state a precise variant of this claim,
taking into account that the generator is based at a particular vertex. 
As before, for a vertex $j\in V$, let 
$f_j\in A$ be the idempotent based at $j$, and $P_j=Af_j$ the 
corresponding projective $A$-module.

\begin{proposition} Assume that the algebra~$A$ is 3-Calabi--Yau. 
Then there is a bijection between the set of finitely generated
left $A$-modules, embeddable into~$P_k$ with finite codimension, 
up to $A$-module isomorphism, and the pairs $(M, m)$ of 
finite-dimensional cyclic $A$-modules $M$ with generator $m\in f_kM$, up to 
$A$-module isomorphism respecting generators.
\label{ideals}
\end{proposition}
\begin{pf} As before, the cyclic module $(M, m)$ defines a surjection 
$\bar m\colon A\rightarrow M$ sending $1\in A$ to $m\in M$. As $m\in f_kM$, 
the kernel $I=\ker\bar m$ decomposes as $I=I_k\oplus \bigoplus_{j\neq k} P_j$, with 
the left $A$-ideal $I_k=If_k$ embedded in $P_k$ with finite codimension.
Isomorphic pairs give rise to isomorphic $A$-modules. Conversely, suppose
that $N\in\AMod$ is embeddable into $P_k$ with finite-dimensional cokernel. 
Take such an embedding $i\colon N\hookrightarrow P_k$, and consider the exact 
sequence of left $A$-modules
\[ 0 \to N \stackrel{i}\to P_k \to M \to 0.\]
Part of the corresponding long exact sequence of abelian groups reads
\[ \Hom_A(M, P_k) \to \Hom_A(P_k, P_k) \to \Hom_A(N,P_k) \to \Ext^1_A(M, P_k).
\] 
However, $M$ is finite dimensional, so by the Calabi--Yau duality, for $i<3$ 
\[\Ext^i_A(M,P_k)\cong \Ext^{3-i}_A(P_k, M)^* = 0,\]
since $P_k$ is a projective $A$-module. Hence
\[ \Hom_A(N,P_k) \cong\Hom_A(P_k, P_k)\cong f_kP_k,\]
where the last isomorphism is given by $g\mapsto g(f_k)\in f_kP_k$.
Thus the only other $A$-maps from $N$ to $P_k$ are multiples of~$i$ by 
some $a\in f_kP_k$. These other embeddings have cokernel
$P_k/Na$ which has $P_k/P_ka$ as quotient, the latter being infinite-dimensional
unless $a\in f_kP_k$ is a constant multiple of the idempotent $f_k$.  
Hence the only embeddings of $N$ into $P_k$ with finite-dimensional quotient
are constant multiples of~$i$, and hence $i$, and thus the surjection
$A\to M$ defining a pair $(M,m)$, are uniquely determined by $N$ up to
$A$-isomorphism.
\end{pf}

As usual, one can presumably promote the statement to an isomorphism of moduli 
spaces, once the relevant moduli problem is formulated for ideals. Since 
this is not relevant for the rest of the paper, I will not pursue this 
direction.

\subsection{Numerical invariants}

Given a scheme $Z$ which admits a symmetric perfect obstruction theory, 
Behrend~\cite{behrend} proves that there is a canonical constructible $\Z$-valued 
function $\nu_Z$ on $Z$. This allows one to define the virtual number of points of 
$Z$ to be the integer
\[ \#^{\rm vir}(Z)  = \sum_{n\in \Z} n \chi(\nu_Z^{-1}(n)),\]
where $\chi$ denotes the topological Euler characteristic.
For the case of the non-commutative algebra $A$, we therefore get integers
\[ D_{k,\v} = \#^{\rm vir}(\M_{i,\v}),\]
depending on a choice of base vertex $k\in V$ and dimension vector 
$\v\in\Z^V$. These can be encoded in the partition function
\[ Z_{A,k}(\q)= \sum_{\v\in\Z^V} D_{k,\v} \q^\v,\]
using a set of auxiliary variables $\q=\{q_i:i\in V\}$. 

\subsection{The commutative Hilbert scheme}\label{sec_commHilbsch}
Recall the example discussed in Example~\ref{ex_threevar}, with 
the quiver~$Q$ having one vertex $V=\{1\}$ and three loop edges, leading to 
$A\cong\C[x_1, x_2, x_3]$ commutative with dimension lattice $\Z$. 
For $n\geq 1$, the moduli space $\M_{1,n}$ is well known to be 
$\Hilb^{[n]}(\C^3)$, the Hilbert scheme of~$n$ points on~$\C^3$. This space
carries a symmetric obstruction theory via its manifestation as the moduli
of rank-1 sheaves with trivial determinant~\cite{th_th, th, bf} on~$\C^3$. 
Theorem~\ref{th_obs} above gives a
new, in some sense more elementary, construction of this obstruction theory.  

The series $Z_{A,1}(q)$ can be computed by torus 
localization~\cite{orv, mnop1, bf}, 
the torus-fixed points being parametrized by 3-dimensional partitions
$\alpha\subset \N^3$. This gives
\[ Z_{A,1}(q) = \sum_\alpha (-q)^{\wt(\alpha)},
\]
where $\wt(\alpha)$ is the total weight of the partition~$\alpha$. 
By a classical result of MacMahon, this function is given by 
\[ Z_{A,1}(q)= M(-q) = \prod_{n\geq 1} (1-(-q)^n)^{-n}. 
\]
Recall also~\cite[Section 5.4]{orv} that finite 3-dimensional partitions 
are in one-to-one
correspondence with dimer configurations in the hexagonal lattice
(honeycomb dimers), with a certain boundary condition.

In the rest of the paper, I study the partition function $Z_A(\q)$ 
for a more complicated example, obtaining analogous results.

\section{The non-commutative conifold}

\subsection{The algebra}\label{sec_conifoldA}

\begin{figure}[ht]
\centering
\ifthenelse{\equal{\pics}{true}}
{\input{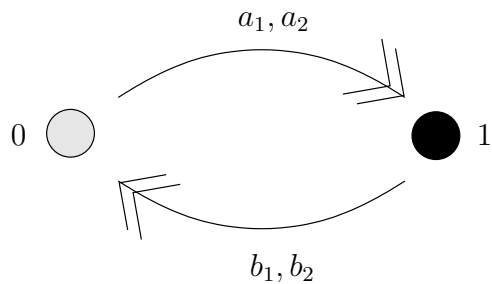}}{}
\caption{The conifold quiver} 
\label{fig_arrows}
\end{figure}

Consider the quiver $Q=\{V,E\}$ of Figure~\ref{fig_arrows}, 
with two vertices $V=\{0, 1\}$, 
four oriented edges $E=\{a_1,a_2\colon 0\to 1, b_1,b_2\colon 1\to 0\}$, and 
relations coming from the quartic superpotential $W=a_1b_1a_2b_2-a_1b_2a_2b_1$
(the Klebanov--Witten superpotential~\cite{kw}). Thus the algebra~$A$
contains the idempotent ring $\C[f_0, f_1]$, and can be given by generators and
relations as
\[ A = \C[f_0, f_1]\langle a_1, a_2, b_1, b_2\rangle / \langle b_1a_ib_{2}-b_2a_ib_1, a_1b_ia_2-a_2b_ia_1 : i=1,2\rangle.
\]
The dimension lattice is~$\Z^V =\Z^2$.
The center $R=Z(A)$ is spanned by $x_1=a_1b_1+b_1a_1, x_2=a_2b_2+b_2a_2, x_3=a_1b_2+b_2a_1, x_4=a_2b_1+b_1a_2$, so 
\[R = \C[x_1,x_2,x_3,x_4]/(x_1x_2-x_3x_4),\]
the ring of functions on the threefold ordinary double point or conifold 
singularity~$Z=\spec(R)$. Indeed, $A$ is a 3-Calabi--Yau algebra, and a crepant 
non-commutative resolution~\cite{vdb2} of the singular variety~$Z$.

\subsection{Torus actions}
\label{sec_toract}

Let $T_F=(\C^*)^E$ be the ``flavour'' torus, acting by diagonally rescaling 
the edge variables. The relations are automatically preserved, and thus $T_F$ 
acts as a group of automorphisms of the algebra~$A$. $T_F$ has two distinguished 
subtori. As a subgroup of $\aut(A)$, it contains the ``gauge'' torus 
$(\C^*)^V/\C^*$, acting on~$A$ by inner automorphisms as 
$r \mapsto (\sum \mu_i^{-1} f_i)\, r \, (\sum \mu_i f_i)$; 
note diagonal elements act trivially. On the other hand, 
$T_F$ also acts on the whole quiver algebra $\C Q$, and hence it has a 
subtorus $T_{F,W}\subset T_F$ which stabilizes the superpotential~$W$. 

Recall the space $X=\{((\phi_e)_{e\in E}, m)\}$ 
of cyclic representations of the quiver~$Q$ on a collection 
of vector spaces $(U_i)_{i \in V}$, generated by a vector $m\in U_k$. 
The flavour torus $T_F=(\C^*)^E=\{(\lambda_e)_{e\in E}\}$ acts on~$X$ by
\[(\lambda_e)\circ((\phi_e), m)= ((\lambda_e\phi_e), m).\]
The torus $T_G=(\C^*)^V=\{(\mu_i)_{i\in V}\}$ also acts on~$X$, 
as a subgroup of $G=\prod_{i\in V}\GL(U_i)$, by 
\[(\mu_i)\circ((\phi_e), m)=((\mu_{h(e)}^{}\mu_{t(e)}^{-1}\phi_e), \mu_k m).\]
The intersection of the image of~$T_F$ and the image of~$G$ 
in $\aut(X)$ is generated by the subtorus 
$T_k=(\C^*)^{V\setminus\{k\}}$ of $T_F$, acting at the 
other vertex $j\in V\setminus\{k\}$.
Hence the moduli space $\M_{k,\v}$ of framed cyclic representations of~$A$
is acted on by the quotient torus $T=T_F/T_k$ of rank~$3$.
The latter has a subtorus $T_W=T_{F,W}/T_k\cap T_{F,W}$ which stabilizes
the regular function $\bar w\in H^0(\O_N)$ used in the construction of
the perfect obstruction theory on $\M_{k,\v}\subset N$.

Explicitly, $T$~can be described as the quotient
of the four-dimensional torus $(\C^*)^4$ acting on the edge variables 
$a_i,b_j$ by the subtorus $\C^*$ acting by weights $(-1, -1, 1,1)$.
$T_W\subset T$ is defined by the condition that the product of the elements
is~$1$.

\subsection{The resolution as moduli space}
Let us construct one well known moduli space of cyclic $A$-modules, that
based at the vertex $0$, corresponding to the dimension vector $(1,1)\in\Z^2$. 
Choosing bases leads to $U_i\cong\C$ and then the arrows $a_1,a_2,b_1,b_2\in\C$ 
as well as the generator $m\in\C^*$ are scalars; cyclicity requires 
that not both $a_1,a_2$ should be zero. Thus
\[\begin{array}{rcl}
\M_{0,(1,1)} & \cong & \left(\C^4_{a_1,a_2,b_1,b_2}\times\C^*_m\setminus\{a_1=a_2=0\}\right)\Big/\C^*(-1, -1, 1,1; 1)\times\C^*(1,1,-1,-1;0)\\
& \cong & \left(\C^4_{a_1,a_2,b_1,b_2}\setminus\{a_1=a_2=0\}\right)\Big/\C^*(-1, -1, 1,1)
\end{array}\]
which is well known to be isomorphic to the crepant resolution~$X$ of the 
conifold~$Z=\spec(R)$. 

\subsection{Pyramid partitions} 
\begin{figure}[ht]
\centering
\ifthenelse{\equal{\pics}{true}}
{\input{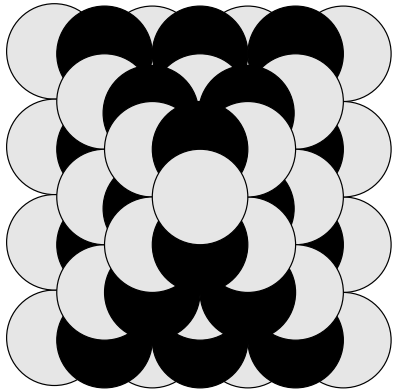}}{}
\caption{The pyramid arrangement} 
\label{fig_pyr}
\end{figure}

Consider the infinite 
combinatorial arrangement~\cite{kenyon} on Figure~\ref{fig_pyr}, 
with two types of layers of stones. 
For $i\geq 0$, there are $(i+1)^2$ stones labelled $0$ and
coloured grey on layer $2i$. On layer $2i+1$, there are $(i+1)(i+2)$ stones 
labelled $1$ and coloured black.

A finite subset $\pi$ of the combinatorial arrangement is a 
{\it pyramid partition}, if for every stone in $\pi$, the (usually two) stones 
immediately above it, of different colour, are also contained in $\pi$. 
For a pyramid partition $\pi$, define $\wt_i(\pi)$ to be the number 
of stones labelled $i$ in $\pi$ for $i=0,1$. 
Denote by $\sP$ the set of all pyramid partitons. 

\begin{remark}\rm The subset of the combinatorial arrangement 
of Figure~\ref{fig_pyr} consisting of grey stones forms
a finitely generated semigroup $\sigma_\Z$, the intersection of $\Z^3$ with 
a polyhedral cone~$\sigma$. The semigroup ring $\C[\sigma_\Z]$ is nothing else but
the ring~$R$, the center of $A$; this is the toric description of the singular
conifold. 

The effect of the non-commutative resolution is to introduce the
odd layers of stones in the picture. The full pyramid defines the 
projective left $A$-module $P_0=Af_0$, where $f_0\in A$ is the idempotent 
corresponding to the chosen vertex $0\in V$. Namely, 
$P_0$ has, as a $\C$-vector space, 
a basis labelled by all the stones in the pyramid, 
and the multiplication by the basic monomials $a_1, a_2, b_1, b_2$
is also encoded (see Figure~\ref{fig_module} below). 
This is analogous to the way the non-negative
octant in the lattice $\Z^3$ gives a basis for the commutative algebra 
$\C[x_1, x_2, x_3]$, with the monomials $x_1, x_2, x_3$ multiplying along 
the edges (compare Example~\ref{ex_threevar} as well as 
\textsection\ref{sec_commHilbsch}).
\end{remark}

\subsection{Torus-fixed points and pyramid partitions}
\label{sec_pyr}
 
The two vertices of the conifold quiver $Q=\{V,E\}$ are symmetric under an
outer automorphism of the algebra; thus I concentrate on cyclic $A$-modules
based on the vertex~$0\in V$, and drop the index~$0$ from the notation.
(See the end of \textsection\ref{sec_br} for more discussion of this choice.)

Consider the action of the rank-two torus $T_W$ of~\textsection\ref{sec_toract} 
on some moduli space $\M_\v$ of framed cyclic $A$-modules. 

\begin{proposition} 
\begin{enumerate}
\item There are finitely many $T_W$-fixed points on the moduli space $\M_{\v}$. 
There is a one-to-one correspondence between fixed points and
pyramid partitions $\pi\in\sP$ of weight $(\wt_1(\pi), \wt_2(\pi))=\v\in\Z^2$. 
\item At each fixed point $P\in \M_{\v}$, the Zariski tangent space to $\M_{\v}$ 
at $P$ has no $T_W$-invariant subspace.
\end{enumerate}
\label{prop_fixedpt}
\end{proposition}
\begin{pf} Suppose that $[(M,m)]\in\M_{\v}$ is an isomorphism class 
of framed cyclic~$A$-modules fixed by the torus~$T_W$, represented by a
framed cyclic $A$-module~$(M, m)$ which can be assumed to be non-zero.
Let $I\lhd A$ be the annihilator of the cyclic generator $m\in M$, 
a left ideal in~$A$. As in~\textsection\ref{sec_ideals}, since the generating 
vector $m\in M$ is based at the vertex $0$, we can write $I=I_0\oplus P_1$, 
where $P_1=Af_1$ consists of all paths starting at the vertex $1$. 
First I claim that~$I_0$ is a monomial ideal: its generators over 
$\C[f_0]$ are monomials in the edge variables $a_i,b_j$. 

To show this, recall that explicitly, as acting on ideals in~$A$, 
the rank-two torus $T_W$ manifests itself as a quotient 
\[T_W=T_{F,W}/\C^*\]
of the rank-three torus
\[ T_{F,W} = \{(\lambda_1, \lambda_2, \mu_1, \mu_2) : \lambda_1\lambda_2\mu_1\mu_2 =1\}\]
acting on the variables $(a_1, a_2, b_1, b_2)$, by the ``gauge'' subtorus
\[\C^*=\{(\lambda,\lambda, \lambda^{-1},\lambda^{-1})\}\subset T_{F,W}.\]  
Generators of $I_0$ split up into sums of paths with the same starting point 
and endpoint, and on elements of~$\C Q$ represented by such paths, this subtorus 
has a constant diagonal action. Hence the covering torus $T_{F,W}$ also acts
on ideals, and the $T_W$-fixed points in $\M_\v$ are in one-to-one correspondence
with $T_{F,W}$-fixed left $A$-ideals $I=I_0\oplus P_1$. 

If the left $A$-ideal $I=I_0\oplus P_1$ is fixed by $T_{F,W}$, then it must be 
generated by $T_{F,W}$-eigenvectors. Let $r(a_1, a_2, b_1, b_2)$ be a (non-commutative)
polynomial in the generators of~$A$,  
which is a $T_{F,W}$-eigenvector in $I_0$. Note that, using the 
relations coming from the superpotential $W$,
\[b_2a_2b_1a_1 = b_2 a_1 b_1 a_2 = b_1a_2b_2a_1=b_1a_1b_2a_2\in A.\]
This element $c=b_ia_jb_{3-i}a_{3-j}\in P_0$ commutes with all elements of $P_0$, and generates 
the weight-$0$ eigenspace of $T_{F,W}$ acting on $P_0$. 
Moving the highest possible power of~$c$ to the right 
in each monomial making up $r$, it follows that
\[r(a_i, b_j)=q(a_1, a_2, b_1, b_2) \cdot p(c)\in A,\] 
where~$p$ is a polynomial with nonzero
constant term, and $q$ is a monomial in the generators $a_i, b_j$. 
Let $J = I_0 \cap Z(A)f_0$, 
where $Z(A)$ is the center of $A$, and $f_0$ is the idempotent at the vertex~$0$. 
Then~$J$ is an ideal in $Z(A)f_0\cong Z(A)$, the coordinate ring of the 
conifold singularity. Since~$I$ is fixed by~$T_W$, the ideal $J$ has to be fixed as well, 
and hence the zero-set of~$J$ must be supported at the singularity
$0\in\spec Z(A)$. Thus this 
zero-set is disjoint from the zero-set of $p(c)\in Z(A)f_0$. By the Nullstellensatz,
this implies that $\langle p, J\rangle = Z(A)f_0$, and hence $q\in I_0$. So
indeed, $I_0$ is generated by monomials. 

Consider the generator~$m\in M$, and the set of vectors 
\[{\mathcal S}=\{m, a_1m, a_2m, b_1a_1m, b_2a_1m, b_1a_2m, b_2a_2m, \ldots\}\subset M.\] 
Since $m$ is a generator, the nonzero vectors in ${\mathcal S}$ form a spanning set
for $M$. Certain sets of vectors in ${\mathcal S}$ will be equal in $M$ as a 
consequence of the relations among the $a_i, b_j$; delete all but one 
vector from each equivalence class. 
Since all remaining relations, contained in the ideal~$I$, are monomial, the 
remaining set of nonzero vectors is linearly independent
and hence forms a basis of~$M$, finite since~$M$ is 
finite dimensional. These nonzero vectors form a finite pyramid 
partition $\pi\in\sP$. 

Conversely, given a pyramid partition $\pi\in\sP$, let $M_i$ for $i=0,1$ be
the $\C$-vector space spanned by basis 
vectors $e_s$ for each stone $s\in\pi$ of label~$i$. The arrows $a_i,b_j$ act by 
mapping between basis vectors as in Figure~\ref{fig_module}. It is immediately seen
that the relations in~$A$ are satisfied. Hence this rule defines an 
$A$-module structure on $M=M_0\oplus M_1$, generated by $m=e_0\in M_0$ corresponding
to the topmost stone of label~$0$. The action of the torus~$T_W$ can be compensated
for by a change of basis, and hence all these modules are in the fixed locus
$(\M_\v)^{T_W}$. This gives the inverse correspondence, concluding the proof of~(1).

\begin{figure}
\centering
\ifthenelse{\equal{\pics}{true}}
{\input{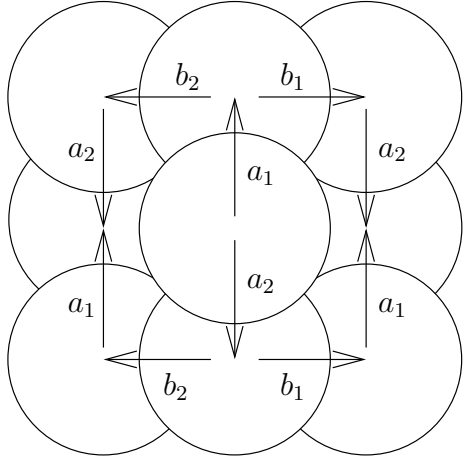}}{}
\caption{A pyramid partition defines a module} 
\label{fig_module}
\end{figure}

To prove (2), by Corollary~\ref{cor_tangent}, it remains to study 
the $T_{F,W}$-action on the space $\Hom_A(I_0, M)$ for a monomial ideal 
$I=I_0\oplus P_1$ with quotient $M=A/I$. I follow the argument 
of~\cite[Lemma 4.1]{bf}. It is enough to show that under the induced action 
of the ``flavour'' torus $T_F = (\C^*)^E$, no weight is a multiple of $(1,1,1,1)$. 

Suppose therefore that $\phi: I_0\to M$ is an eigenvector with weight 
$w(1,1,1,1)$. Suppose first that $w\geq 0$. 
As the $(1,1,1,1)$ eigenspace of $P_0$
is spanned by the element $c=b_ia_jb_{3-i}a_{3-j}$ already used above,
for any element $a\in I_0$, necessarily
\[\phi(a)\equiv c^w a \equiv 0 \mod I,\] and so $\phi=0$.

To treat the case $w<0$, let $\alpha$ be the smallest positive integer such that 
$(b_1a_1)^\alpha\in I_0$, and let $\beta$ be the smallest positive integer such that
$(b_2a_2)^\beta(b_1a_1)^{\alpha-1}\in I_0$. Since $\phi$ is an $T_W$-eigenvector with 
weight $w(1,1,1,1)$,  
\[\phi((b_2a_2)^\beta(b_1a_1)^{\alpha-1})\equiv(b_2a_2)^{\beta+w}(b_1a_1)^{\alpha-1+w} \mod I.
\]
Now compute, using the fact that $b_1a_1$ commutes with $b_2a_2$:
\begin{eqnarray*}\phi((b_2a_2)^\beta(b_1a_1)^\alpha)=b_1a_1\phi((b_2a_2)^\beta(b_1a_1)^{\alpha-1})& \equiv &  b_1a_1(b_2a_2)^{\beta+w}(b_1a_1)^{\alpha+w-1} \\
& \equiv & (b_2a_2)^{\beta+w}(b_1a_1)^{\alpha+w} \mod I.
\end{eqnarray*}
On the other hand,
\[\phi((b_1a_1)^{\alpha})\equiv 0 \mod I,
\]
since there is no monomial in~$P_0$ with negative $T_W$-weights. Hence 
\[\phi((b_2a_2)^\beta(b_1a_1)^\alpha)=(b_2a_2)^\beta\phi((b_1a_1)^\alpha) \equiv 0 \mod I.
\]
Comparing the two expressions, 
\[(b_2a_2)^{\beta+w}(b_1a_1)^{\alpha+w} \in I. 
\]
Since $w<0$, this contradicts the definition of $\beta$. 
\end{pf}

I will later need the following additional information on fixed points.
Recall from~\textsection\ref{sec_modconst} 
the construction of the moduli space $\M_{\v}$ 
of framed cyclic modules as the free quotient of an 
locally closed subset $X\subset S^0\subset S$ of a vector space $S$ of linear maps 
by a group $G$. 

\begin{lemma} Let $\pi\in\sP$ be a pyramid partition, and $M_\pi$ the
framed cyclic module defined by~$\pi$. Then at the point 
$[M_\pi]\in X$, the parities of the tangent spaces $\dim T_{[M_\pi]} X$ and 
$\dim T_{[M_\pi]}S$ coincide.
\label{lem_par}
\end{lemma}
\begin{pf} 
In the prescription given in the proof of Proposition~\ref{prop_fixedpt},  
the $A$-module $M_\pi$ has a vector space 
basis indexed by the stones of the partition $\pi$. The affine space $S$ inherits
natural coordinates given by the corresponding matrix entries. In this coordinate
system, $[M_\pi]\in S$ has coordinates $0$ and $1$ as dictated by the 
partition $\pi$; in Figure~\ref{fig_module}, the arrows drawn correspond to 
coordinates equal to 1, the rest being 0. 

The embedding $X\subset S$ is locally defined by the superpotential equations 
$d\Tr(W)$. Writing these in the matrix entries, we get
cubic equations for every pair of basis vectors indexed by 
stones $i,l\in\pi$ of opposite colours, and appropriate composable arrows
$a, b, c$ of the conifold quiver, with $a$ and $c$ different, of the form 
\begin{equation}\label{eq_rels}[il; abc]\colon \ \ \ \sum_{j,k} a_{ij}b_{jk}c_{kl} = \sum_{j,k} c_{ij}b_{jk}a_{kl}.\end{equation}
Here $j,k\in \pi$ run over all stones of appropriate colour. To get the 
embedding of tangent spaces, these
equations have to be linearized 
near the point $[M_\pi]\in X$ with all coordinates either~$0$ or~$1$. 
By changing variables to new variables $a_{ij}'=a_{ij}-1$ whenever $a_{ij}=1$,
the linearized form of~(\ref{eq_rels}) is non-zero only 
if at least two of the consecutive arrows $a_{ij}, b_{jk}, c_{kl}$, 
or $c_{ij}, b_{jk}, a_{kl}$ are non-zero. If for some $i$ and $l$,
there are three nonzero arrows $a_{ij}, b_{jk}, c_{kl}$, then (by the relations) 
there is another set $c_{im}, b_{mn}, a_{nl}$ of nonzero arrows, and the 
linearized form is
\[[il; abc]_{\rm lin}\colon \ \ \ a'_{ij} + b'_{jk} + c'_{kl} + \ldots = c'_{im} + b'_{mn} + a'_{nl} + \ldots; \]
If only two of the arrows, say $a_{ij}$ and $b_{jk}$ are nonzero, the
linearized form is
\[ [il; abc]_{\rm lin}\colon \ \ \ c_{kl} + \ldots = 0.\]
In both cases, $\ldots$ represents further possible terms of the same shape.  

\begin{figure}
\centering
\ifthenelse{\equal{\pics}{true}}
{\input{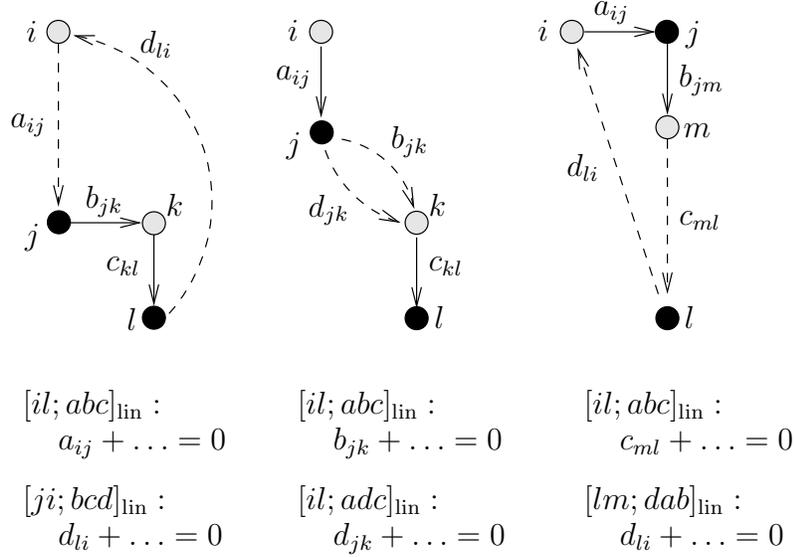}}{}
\caption{Configurations of arrows, and paired linear equations}
\label{fig_confs}
\end{figure}

I now claim that it is possible to pair these linearized equations so that 
the embedding $T_{[M_\pi]} X\hookrightarrow T_{[M_\pi]}S$ is cut 
out by an even number of linearly independent equations, proving the lemma. 
The proof proceeds by induction on the size of the pyramid partition~$\pi$, 
the claim being obvious for small partitions. 

Take a stone $l\in\pi$ so that $\pi\setminus\{l\}$ is also
a partition. Assume first that the stone $l$ only touches one other stone 
of $\pi$. There is then only one nonzero incoming arrow $c_{kl}$ to $l$, 
and $l$ is not at the end of a chain of three consecutive nonzero arrows
$a_{ij}, b_{jk}, c_{kl}$, since (by the relations) that would necessitate
another nonzero incoming arrow $a_{nl}$. It is now easy to check that there 
are only three possible configurations of two nonzero arrows,
leading to new nonzero linear relations, being part of a path involving $l$. 
The three configurations are depicted on Figure~\ref{fig_confs}, together with
the pairs of linear equations that the configurations give rise to. 
In the figure, solid arrows represent 1's, whereas dotted arrows represent 0's.
For example, in the middle case, we have two different nonzero
arrows $a_{ij}, c_{kl}$, and there are two ways of completing the pair to a path: 
$a_{ij}b_{jk}c_{kl}$ and $a_{ij}d_{jk}c_{kl}$, leading to two linear relations
cutting out the linear subspace of interest.

It is immediate that in all cases the linear equations paired are indeed linearly
independent, and the pairs are well-defined. This proves the claim, and hence the
lemma, in this case. 
The second case, where there are two nonzero incoming arrows $c_{kl}$ and $a_{nl}$
to~$l$, is similar, with a slightly larger number of diagrams. The details are left
to the reader.
\end{pf}

\begin{corollary} Let $\v=(d_0, d_1)$ be a dimension vector. 
The parity of the dimension of the Zariski tangent space
at a $T_W$-fixed point $P\in(\M_{\v})^{T_W}$ is the same as the parity of $d_1$.
\label{cor_dim}
\end{corollary}
\begin{pf} Let $Q\in X$ lie above $P\in\M_{\v}$. Since the $G$-action is free, 
\[ \dim T_P\M_{\v} = \dim T_Q X - \dim G.\]
On the other hand, using Lemma~\ref{lem_par} and the fact 
that $S$ is just a vector space, 
\[ \dim T_Q X \equiv \dim T_Q S  = \dim S \ (\mod\ 2);\]
hence
\[ \dim T_P\M_{\v} \equiv \dim S - \dim G \ (\mod\ 2).\]
Looking at the definitions of $S, G$, this difference of dimensions is indeed
\[ 4 d_0 d_1 + d_0 - d_0^2 - d_1^2 \equiv d_1 \ (\mod\ 2).\]
\end{pf}

\subsection{Pyramid partitions and dimer configurations}
\label{sec_dim}

This section is not necessary for the logical flow of the paper. However, 
it points the way to generalizations, so I include it here. 

Let $\L\cong\Z^2$ be the square lattice, thought of as an unoriented 4-regular 
bipartite graph with infinite vertex set $\V$ and edge set $\E$. Recall that 
a dimer configuration on the square lattice is a subset $\F\subset\E$ of 
the edges, so that every vertex $v\in\V$ is incident to exactly
one edge in $\F$. Given a fixed dimer configuration $\F_1$, a dimer configuration
is said to be asymptotic to $\F_1$, if the set $\F\setminus \F_1$ is finite.

\begin{figure}
\centering
\ifthenelse{\equal{\pics}{true}}
{\input{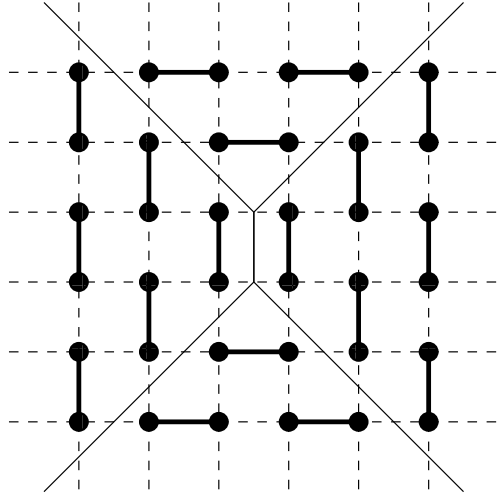}}{}
\caption{The ``length one empty room'' square dimer configuration $\F_1$}
\label{fig_empty}
\end{figure}

Let $\F_1$ be the ``length one empty room'' square dimer configuration shown 
on Figure~\ref{fig_empty}. 

\begin{figure}
\centering
\ifthenelse{\equal{\pics}{true}}
{\input{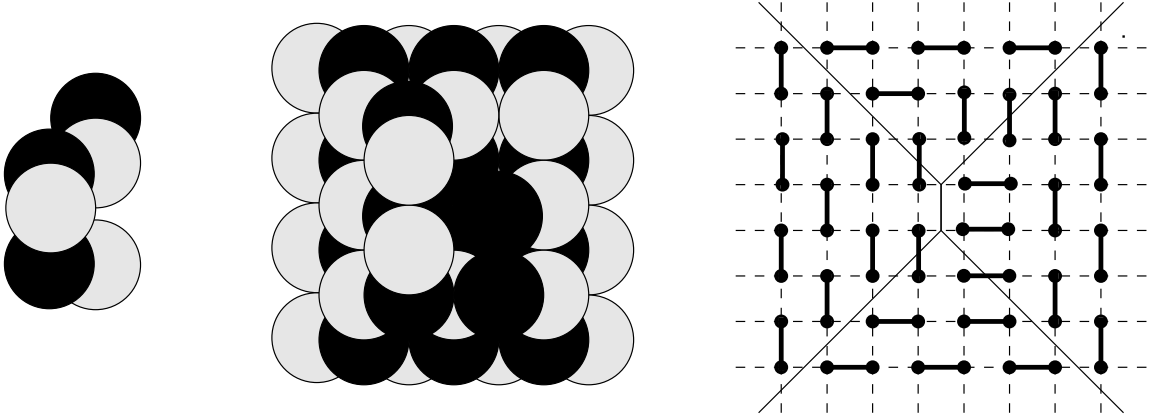}}{}
\caption{A pyramid partition, its complement with the grid, and the resulting dimer 
configuration}
\label{fig_corresp}
\end{figure}

\begin{proposition} There is a one-to-one correspondence between 
pyramid partitions $\alpha\in\sP$ and dimer configurations $\F$ 
on the square lattice $\L$ asymptotic to the configuration $\F_1$.
\label{prop_dimers}
\end{proposition}
\begin{pf} 

This is analogous to the correspondence between honeycomb dimers and
ordinary 3-dimensional partitions. Given a pyramid partiton $\alpha\in\sP$, consider
its complement in the infinite pyramid arrangement. Now associate to this complement
a dimer configuration a~$\F$ as follows. Looking from above, impose a square grid
in which the squares lie over the balls as on Figure~\ref{fig_corresp}. Now consider 
those edges which lie over balls or half-balls visible from above. As illustrated, 
black balls or half-balls correspond to horizontal dimers, whereas grey ones 
correspond to vertical dimers. The empty pyramid partition gives the ``length
one empty room'' $\F_1$. Since $\alpha\in\sP$ is finite, $\F$ is asymptotic to
$\F_1$ and it is easy to see that all such arise.
\end{pf}

Different pictures, perhaps more illuminating for some, can be found 
in~\cite[Figures 1 and 2]{young}.

\begin{remark}\rm 
A general periodic dimer model, under some extra conditions, also gives rise 
to a non-commutative toric algebra~$A$ defined by a superpotential. 
Moreover, there is a correspondence between torus-fixed points in moduli 
spaces $\M_{k,\v}$ and dimer configurations with fixed asymptotic behaviour. 
This will be discussed in future work.
\end{remark}

\subsection{The non-commutative partition function}
\label{sec_prod}
Consider the generating function 
\[ Z_{A}(\q)= \sum_{\v\in\Z^2} D_{\v} \q^\v\]
of virtual counts of framed cyclic representations of $A$. 

\begin{theorem} The generating function $Z_{A}(\q)$ can be expressed 
combinatorially as
\[
Z_A(q_0, q_1) = \sum_{\pi\in\sP} (-1)_{}^{\wt_1(\pi)} q_0^{\wt_0(\pi)}q_1^{\wt_1(\pi)},
\]
where $\sP$ is the set of finite pyramid partitions, or equivalently 
by Proposition~\ref{prop_dimers}, the set of square dimer configurations~$\F$
asymptotic to the ``length one empty room''~$\F_1$.
\end{theorem}
\begin{pf} The rank-2 torus $T_W$ acts on all
moduli spaces~$\M_{\v}$. By construction, the symmetric obstruction theory 
constructed in Theorem~\ref{th_obs} is $T_W$-equivariant in the sense of~\cite{bf}. 
Since $D_{\v}$ is a weighted Euler characteristic, with  
weights which are constant on $T_W$-orbits, nontrivial orbits do
not contribute, and hence $D_{\v}$ receives contributions from torus-fixed
points only. At a $T_W$-fixed point $P\in (\M_\v)^{T_W}$, the contribution is simply 
$(-1)^{\dim T_P\M_\v}$ by applying~\cite[Theorem 3.4]{bf} to a sufficiently 
general one-dimensional subtorus of~$T_W$. Hence, the statement follows from
Corollary~\ref{cor_dim}.
\end{pf}

The following result was conjectured in an earlier version of this paper,
on the basis of extensive computational evidence, generalizing the one-variable 
specialization conjectured earlier by Kenyon~\cite{kenyon}.

\begin{theorem} \label{thm_young}
{\rm (Young~\cite{young})} The partition function $Z_A(\q)$ 
admits the following infinite product expansion.
\begin{equation}
\label{fnA}
Z_A(q_0, q_1) = M(-q_0q_1)^2 \prod_{k\geq 1}\left(1+q_0^k(-q_1)^{k-1}\right)^k\left(1+q_0^k(-q_1)^{k+1}\right)^k.
\end{equation}
\end{theorem}

\section{Interpretations and generalizations}
\label{chap_interp}

\subsection{The commutative partition function}
\label{sec_cdt}

Recall the rank one Donaldson--Thomas partition function~\cite{mnop1}
of the smooth commutative Calabi--Yau resolution $X$, the 
resolved conifold. The resolution $\pi\colon X\to Z$ contracts
a single rational curve~$C\subset X$ to the singular point of~$Z$.
Let $D_{n,d}$ be the DT invariant~\cite{th}
of ideal sheaves~$\I$ on~$X$ with $\chi(\I)=n$ 
and homology class $d[C]$. Define the DT partition function 
of~$X$ as the series
\[Z_X(q,t) =\sum_{n\in\Z}\sum_{d\geq 0} D^n_d q^n e^{-dt}.\]

\begin{figure}
\centering
\ifthenelse{\equal{\pics}{true}}
{\input{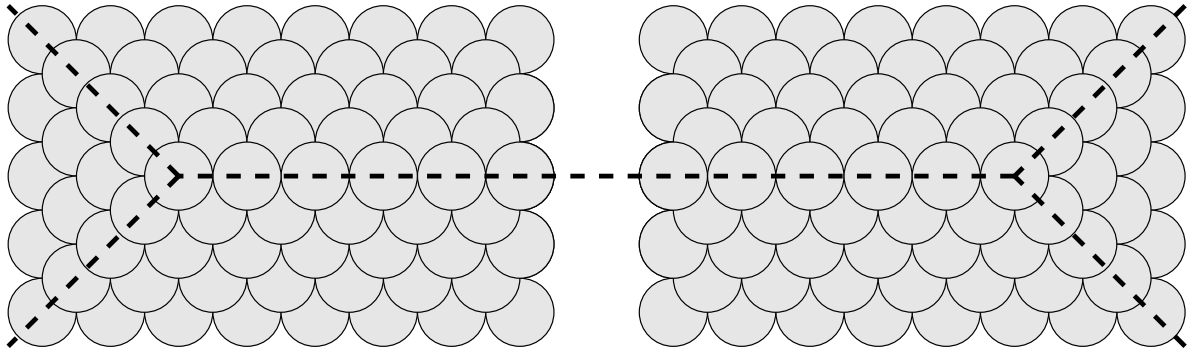}}{}
\caption{The combinatorial arrangement for the commutative resolution, with the
toric web diagram indicated}
\label{fig_comm}
\end{figure}

An argument using torus localization, similar to the one given above,
gives the following combinatorial interpretation of these DT invariants.
Up to sign, $D^n_d$ is the number of pairs of semi-infinite 3-dimensional 
partitions $(\alpha_1,\alpha_2)$, bounded in the direction 
of two coordinate axes, with infinitely long necks in the third direction with 
a common ordinary partition $\lambda\vdash d$ as cross-section, and having total 
(renormalized) volume $|\alpha_1| + |\alpha_2|=n-f(\lambda)$, where $f(\lambda)$ is
a certain combinatorial invariant~\cite[Lemma 5]{mnop1} of the partition~$\lambda$.
In other words, $D^n_d$ counts
certain subsets of the arrangement of Figure~\ref{fig_comm}; note that the shape of 
the latter is closely related to the toric combinatorics of the resolved 
conifold~$X$ (indicated by the dashed line). Consequently~\cite{orv, inov, mnop1}, 
one obtains
\begin{theorem} 
\begin{equation}
\label{fnX}
Z_X(q,t) = M(-q)^2 \prod_{k\geq 1} \left(1-(-q)^k e^{-t}\right)^k.
\end{equation}
\end{theorem}

In fact, the singularity $Z=\spec(R)$ admits two, isomorphic, crepant 
resolutions $X, X^+$ related by a flop $X^+\dashrightarrow X$. Under the 
natural isomorphism $H^2(X,\R)\to H^2(X^+, \R)$ 
induced by the flop, the positive classes have opposite sign.  
Therefore the DT partition function of $X^+$ can be written in the variables
$q, t$ as
\begin{equation}
\label{fnX+} 
Z_{X^+}(q,t) = M(-q)^2 \prod_{k\geq 1} \left(1-(-q)^k e^{t}\right)^k.
\end{equation}

\subsection{Partition functions as functions on Bridgeland space and wall crossing}
\label{sec_br}

The formulae~(\ref{fnA}), (\ref{fnX}) and~(\ref{fnX+}) are closely related. 
The change of variables $q=q_0q_1$, $z=q_1=e^{-t}$ gives
\begin{eqnarray}\label{eq_ZAprodform} Z_X(q, z) & = & M(-q)^2\prod_{k\geq 1} (1-(-q)^kz)^k,\nonumber \\
Z_A(q, z) & = & M(-q)^2\prod_{k\geq 1} (1-(-q)^kz)^k(1-(-q)^kz^{-1})^k,\\ 
Z_{X^+}(q, z) & = & M(-q)^2\prod_{k\geq 1} (1-(-q)^kz^{-1})^k.\nonumber
\end{eqnarray}
Hence dividing by the MacMahon factors, the reduced partition functions satisfy
the curious formal factorization property
\begin{equation} Z'_A(q,z) = Z'_X(q,z) Z'_{X^+}(q,z).
\label{eq_factor}
\end{equation}

The variable $z=e^{-t}$ coordinatizes a parameter space
naturally associated to the problem. To discuss this space, I need some
definitions. 
Recall the map $\pi\colon X\to Z$, contracting the rational curve~$C$. Let 
$\D(X/Z)$ denote the bounded derived category of coherent sheaves 
on~$X$, supported on a neighbourhood of the exceptional curve~$C$. 
It is known~\cite{vdb} that $\D(X/Z)$ is equivalent to the derived category 
$\D_\nilp(\AMod)$ of complexes of finitely generated $A$-modules with 
locally nilpotent cohomology modules.

Let $\Stab(X/Z)$ denote the component of the space of normalized stability 
conditions~\cite{br, br_ICM} 
on the category $\D(X/Z)$, which contains stability conditions
whose heart is the category of sheaves on $X$ supported along $C$. 
By~\cite{br}, $\Stab(X/Z)$ is a one-dimensional complex manifold. 
Let finally $\DAut(X/Z)$ be the subgroup of the derived autoequivalence 
group of $\D(X/Z)$ which preserves this component. Then it is 
proved in~\cite{br_ICM, toda} that
\[ \Stab(X/Z)/\DAut(X/Z) \cong \PP^1\setminus\{0,1,\infty\}.
\]
In physics language, this is the full complexified K\"ahler structure moduli 
space of~$X$. 

\begin{figure}
\centering
\ifthenelse{\equal{\pics}{true}}
{\input{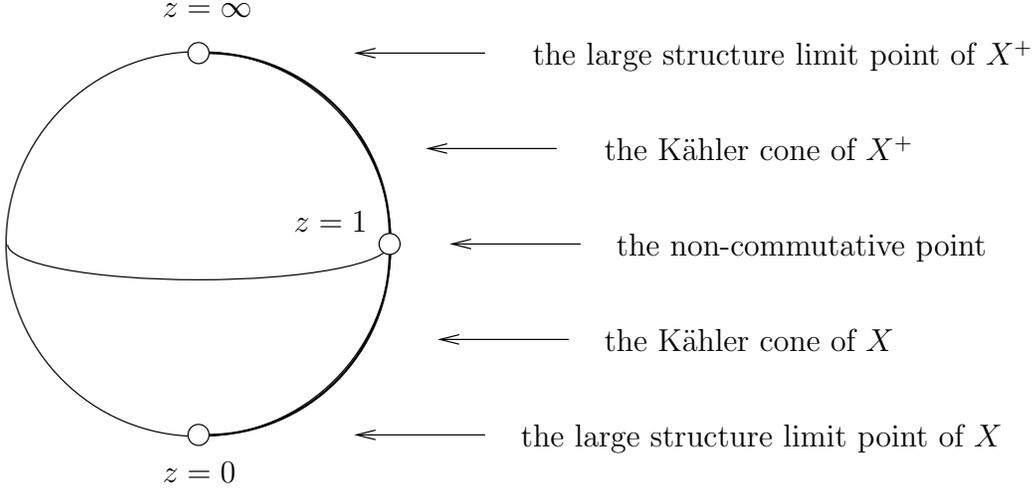}}{}
\caption{The K\"ahler structure moduli space of $X$} 
\label{fig_modspace}
\end{figure}

The space $\Stab(X/Z)/\DAut(X/Z)$
and its structures are pictured on Figure~\ref{fig_modspace}. 
The parameter $t$ naturally lives in the upper or lower half planes 
$\{B + i\omega | \pm\omega>0\}$ of the complexified
cohomology space $H^2(X,\C)\cong\C$, with $\omega=-i\im(t)$ measuring the 
volume of the embedded projective line in~$X$. After taking the quotient 
by $\DAut(X/Z)$, $e^{-t}$ lives in the punctured upper and lower hemispheres 
in $\PP^1\setminus\{0,1,\infty\}$. The two cusps $z=0, \infty$ correspond to the 
``large volume'' limits $t\to\pm\infty$ of $X, X^+$. For 
a stability condition corresponding to a point on the upper or lower hemispheres 
$\{|z|>1\}$, $\{|z|<1\}$, the $t$-structure 
on $\D(X/Z)$ has heart $\Coh(X/Z)$, respectively $\Coh(X^+/Z)$ (up to shift). 
Along the punctured equator $\{|z|=1\}\setminus \{z=1\}$, 
all stability conditions correspond to the perverse $t$-structure on~$\D(X/Z)$,
with heart $\Anilp$, the category of nilpotent $A$-modules. 
In this sense, the cusp $z=1$ is naturally 
associated to the non-commutative algebra~$A$. 

At this point, one would like to say that the functions 
$Z_A$, $Z_X$, $Z_{X^+}$ count
objects from $\D(X/Z)$, which are 
``stable in the limit''. This is however not correct as it stands. 
The commutative DT partition function~$Z_X$ of~$X$ 
certainly counts objects of $\Coh(X/Z)$, at least 
if one thinks of it as counting structure sheaves rather than ideal sheaves. 
However, most of these are not stable: for example,
the sheaf~$E= \O_X/\I_C^2$ on~$C$, where $\I_C$ is the ideal defining 
the curve~$C$, sits in a sequence
\[0\to \O_C(1)^{\oplus 2} \to E \to \O_C \to 0\]
which is a destabilizing subsequence for any stability condition on
$\D(X/Z)$ with heart~$\Coh(X/Z)$.
Similarly, $Z_A$~does count objects of $\Anilp$, since all $T$-fixed
$A$-representations are in fact nilpotent, but most of these are non-stable
(though semistable) in $\D_\nilp(\Amod)$. 

An alternative is to return to the original interpretation of 
the DT partition function~$Z_X$ of~$X$ as a generating function
counting ideal sheaves. These
are indeed Gieseker stable on $X$, but they are not objects of $\D(X/Z)$, only 
of~$\D(X)$. However, very little is known about the space of (normalized) 
stability conditions on~$\D(X)$, though presumably it is also a cover of the
thrice punctured Riemann sphere. It is conceivable that near the large volume
limit point corresponding to~$X$, every Gieseker stable sheaf on~$X$ becomes
Bridgeland stable, so $Z_X$ does in fact count stable objects near a cusp 
in the space of stability conditions. Similarly, by Proposition~\ref{ideals},
the function $Z_A$ can be thought of as counting certain finite-codimension ideals 
in~$A$, and one might hope that in fact these, and only these, objects of
$\D(X)\cong\D(\AMod)$ are stable under an appropriate limiting stability condition
as one approaches the conifold point $z=1$ along the equator. 
Note that there are two ways to do that, perhaps corresponding
to the two choices of vertex for the cyclic generator. 
(Compare~\cite[Section 6.1.2]{dm} for a similar example of ``opposite'' sets of 
objects becoming stable along two opposite rays leading to the same limit point). 
The change from 
$Z_X$ to $Z_A$ to $Z_{X^+}$ might then be interpreted as an instance of 
wall crossing in the space of stability conditions, with a countably infinite 
number of walls between $z=0$ and $z=1$, respectively $z=1$ and $z=\infty$. 
It would be very interesting to pursue this direction further.

Note finally that although $Z_A$ has a Taylor series expansion in the
variables $(q0,q1)=(q/z,z)$, it is only a Laurent series in the variables
$(q,z)$ near $z=0$ or $z=\infty$ with infinitely many positive and
negative powers of $z$. In this sense, this change of variables is
~Snon-perturbative~T. As in~(\ref{eq_factor}), $Z_X$ and $Z_{X^+}$ are the
positive and negative parts of $Z_A$ under a factorization on~$\PP^1_z$,
reminiscent of Birkhoff factorization; these can then be expanded as
Taylor series near the appropriate cusps.

\subsection{Coefficients at small volume and ``modularity''}
\label{sec_mod}

Starting with the form~(\ref{eq_ZAprodform}) of the partition function $Z_A$,
under the substitution $z=e^{-t}$ discussed in the last section, write
\begin{eqnarray}\label{eq_1} 
\log Z_A(q, t) & = & \sum_{k\geq 1} k\cdot \left(\log(1-(-q)^k e^{-t} ) + \log(1-(-q)^ke^t) - 2\log(1-(-q)^k)\right) \nonumber\\
& = & \sum_{k,n \geq 1} \frac{k}{n}\,  (-q)^{kn} \left(e^{-nt} + e^{nt} -2 \right)\nonumber \\
& = & \sum_{d\geq 1} d \cdot (-q)^d \cdot \sum_{n|d} \frac{1}{n^2}\left( 2 \sinh \frac{nt}{2}\right)^2.
\end{eqnarray}
This is formally very similar to the way the change of variables $q=-e^{iu}$ 
of~\cite{mnop1} leads to the Gopakumar--Vafa form~\cite{gv} of the closed string 
partition function. However, it is very different in substance, since~$t$ is the 
K\"ahler parameter. Thus~(\ref{eq_1}) is valid near small volume, or equivalently 
high spacetime curvature, and large string coupling~$u$. Note that it is 
essential to start with the whole~$Z_A$, including its positive and negative parts 
as well as the MacMahon factor, to get this form. 

To continue, use the series expansion 
\[ \left(2\sinh\frac{x}{2}\right)^2 = \sum_{k\geq 1} \frac{2}{(2k)!} \, x^{2k}
\] 
to get
\begin{eqnarray}\label{eq_2} 
\log Z_A(q, t) & = & \sum_{d, k\geq 1} \frac{2d}{(2k)!} \,(-q)^d \,t^{2k} \cdot \sum_{n|d} n^{2k-2} \nonumber \\
& = &  \sum_{d, k\geq 1} \frac{2(-1)^dd \,\sigma_{2k-2}(d)}{(2k)!} \,q^d\, t^{2k},
\end{eqnarray}
with the divisor power sum function for positive integers~$d$ defined by
\[\sigma_s(d) = \sum_{n|d} n^s.
\]
The expression (\ref{eq_2}) has the form of a generating 
series for a set of new ``(connected) enumerative invariants''
\[ M_{d,k} =  \frac{2(-1)^dd\cdot\sigma_{2k-2}(d)}{(2k)!}, \ \ d,k\geq 1,
\] 
associated to~$X$, at small volume/high curvature and large string
coupling. I have no geometric or physical interpretation for 
these numbers at present.

To continue, for integers $l\geq 1$, introduce the series 
\[ E_l(x) = \sum_{d\geq 1} \sigma_{l-1}(d)x^d.
\]
The reason for the notation will become apparent presently. 
(\ref{eq_2}) becomes
\begin{equation}\label{eq_E} 
\log Z_A(q, t) = \sum_{k\geq 1} \frac{2}{(2k)!} \,t^{2k} \, q\frac{d}{dq}E_{2k-1}(-q).\end{equation}
On the other hand, the standard identity (valid for $\im\tau>0$, $l>1$) 
\[  \sum_{n\in \Z}\frac{1}{(\tau+n)^l} = \frac{(-2\pi i)^l}{(l-1)!}\sum_{m\geq 1} m^{l-1} e^{2\pi i m \tau},
\]
leads for $l\geq 2$ to
\[ E_l(e^{2\pi i \tau}) = \frac{(l-1)!}{(-2\pi i)^l} \sum_{m\geq 1, n\in\Z} \frac{1}{(m\tau+n)^l}.
\]
Thus, $E_l$ is ``essentially'' an Eisenstein series. Indeed, 
for $l=2k$ even, we can double the sum to run over $m\in\Z\setminus\{0\}, n\in\Z$ 
and add the appropriate constant term, to get the standard modular (quasi-)invariant
Eisenstein series $G_{2k}$. 
However, for the odd values of interest in~(\ref{eq_E}), the 
signs work against us. So $E_{2k-1}$ is not modular, but it is intriguingly 
close to being so. Compare~\cite[Appendix E]{dmp} for a similar discussion 
involving the MacMahon function.

\subsection{A generalized partition function}\label{sec_higher}

The combinatorial arragement of Figure~\ref{fig_pyr}, leading to pyramid partitions,
has an obvious generalization to an elongated rectangular pyramid with~$n$ 
black stones on level~$0$. In the dimer model, finite subsets of this
arrangement correspond to configurations with the ``length $n$ empty room'' 
asymptotics.

Let $\sP^{(n)}$ denote the set of finite partition-like subsets 
of the elongated pyramid configuration of length $n$; as before, let
$\wt_i(\pi)$ denote the number of stones in $\alpha\in\sP^{(n)}$ of 
colour $i$. Define the partition function
\[ Z_A^{(n)}(q_0, q_1) = \sum_{\pi\in\sP^{(n)}} q_0^{\wt_0(\pi)}(-q_1)_{\ }^{\wt_1(\pi)}.
\]
The choice of signs is motivated by Corollary~\ref{cor_dim},
but this is only by analogy, since the precise enumerative interpretation of 
these invariants is unclear. 
The following generalization of Theorem~\ref{thm_young} was again conjectured
in an earlier version of this paper (compare Kenyon~\cite{kenyon} for the 
one-variable specialization), based on computational evidence:

\begin{theorem} {\rm (Young~\cite{young})}
The partition function $Z_A^{(n)}(\q)$ admits the following infinite product 
expansion.
\begin{equation}
Z_A^{(n)}(q_0, q_1) = M(-q_0q_1)^2 \prod_{k\geq 1}\left(1+q_0^k (-q_1)^{k-1}\right)^{k+n-1}\left(1+q_0^k (-q_1)^{k+1}\right)^{\max(k-n+1,0)}.
\end{equation}
\end{theorem}
\noindent 

In the recent physics literature on the conifold, there are expressions which 
are remarkably similar to the terms appearing as the ratio
$Z_{A}^{(n)}/Z^{\ }_{A}$; compare~\cite[(3.30)-(3.32)]{hy} as well 
as~\cite[(3.6)]{kp}. These computations are in the context of Lagrangian branes, 
so not immediately applicable, but they suggest a relationship between these 
generalized partition functions and higher-rank Donaldson--Thomas theory. 
I~hope to return to this point in future work.

\subsection{The global case}
\label{sec_glob}
In this section, I indicate a possible extension of the ideas of the paper to 
a global context. 
Assume that~$Y$ is a projective variety, singular at a finite set of conifold
points (nodes). Assume further that there is a Calabi--Yau small resolution 
$\pi\colon X \to Y$, together with a line bundle $\L\in\pic(X)$ such that 
$\L|_{C_i}\cong \O_{C_i}(1)$ on all $\pi$-exceptional curves $C_i\cong\PP^1$. 

\begin{proposition} {\rm (Van den Bergh)}
There exists a sheaf $\sA$ of associative, non-commutative 
algebras on $Y$, such that over $Y\setminus{\rm Sing}(Y)$,
$\sA_{Y\setminus{\rm Sing}(Y)}$ is a sheaf of Azumaya algebras (locally,
matrix algebras over the structure sheaf), and in a neighbourhood of each of 
the nodes, sections of $\sA$ define the non-commutative resolution~$A$ discussed 
before. Further, there is a derived equivalence $\D(Y)\cong\D(\sA\rm{-mod})$. 
\end{proposition}
\begin{pf} Define the sheaf $\sA$ of $\O_Y$-algebras by 
\[ \sA = \pi_*\End_X(\O_X\oplus\L). 
\]
This is clearly a sheaf of Azumaya algebras on $Y\setminus{\rm Sing}(X)$, whereas 
by~\cite{vdb}, in a neighbourhood of each of the nodes, sections of $\sA$ define
the noncommutative crepant resolution~$A$. 
The derived equivalence follows from~\cite[Propositions 3.3.1-3.3.2]{vdb}. 
\end{pf}
Thus, $\D(\sA\rm{-mod})$ is a 3-Calabi--Yau category, with Serre functor the third
power of the shift functor, and thus appropriately
defined moduli spaces of (stable) $\sA$-modules will carry a perfect obstruction 
theory. It would be interesting to compare the invariants thus obtained to the 
ordinary DT invariants of~$Y$.

\subsection{Some quotient examples}
\label{sec_orb}
Phenomena similar to the above have been observed by~\cite{by} for some 
orbifolds $\C^3/\Gamma$, 
where the singularity admits a semismall resolution (no divisors over points, 
equivalently all exceptional fibres are rational curves), such as 
$\Gamma<\SL(2,\C)<\SL(3,\C)$ and $\Gamma=\Z/2\times\Z/2<\SL(3,\C)$.
Let $X\to \C^3/\Gamma$ be the canonical crepant resolution given 
by the $\Gamma$-Hilbert scheme. Consider also the canonical non-commutative 
resolution given by the cross product algebra
\[ A = \C[x_1, x_2, x_3] \star \C\Gamma.\]
This algebra is known to be (Morita equivalent to) the quiver algebra of the McKay 
quiver of~$\Gamma$, with a specific superpotential~\cite[Theorem 4.4.6]{ginzburg}.

In complete analogy with the above discussion, there is a commutative 
Donaldson--Thomas partition function~$Z_X$, as well as a noncommutative 
partition function~$Z_A$. Assuming further that $\Gamma$ is abelian, 
both can again be computed by torus localization; the computation of~$Z_A$ 
localizes to a set indexed by standard 3-dimensional partitions 
coloured with characters of the group~$\Gamma$. 

The resulting formulae 
depend on a set of roots $\Sigma\subset\Z^r$, together with a subset 
$\Sigma^+\subset\Sigma$ of positive roots. For a set of dual variables 
${\bf t}$, the commutative partition function takes the form
\[ Z_X(q, {\bf t}) = M(-q)^\chi \prod_{k\geq 1} \prod_{\alpha\in\Sigma^+}(1-(-q)^ke^{-\langle \alpha, {\bf t}\rangle})^{\pm k}
\]
(product over positive roots), whereas $Z_A$ is of the form
\[Z_A(q, {\bf t}) = M(-q)^\chi \prod_{k\geq 1} \prod_{\alpha\in\Sigma}(1-(-q)^ke^{-\langle \alpha, {\bf t}\rangle})^{\pm k}
\]
(product over all roots). In further analogy with the case of the resolved
conifold, the resolution $X$ admits various flops~$X^i$ (geometric or derived) 
which have partition functions~$Z_{X^i}$ obtained by choosing a different 
Weyl chamber as dictated by the flop. In particular, there is a
flop~$X^+$ of~$X$ corresponding to the opposite chamber, together
with a factorization
\[ Z'_A(q, {\bf t}) = Z'_X(q, {\bf t})Z'_{X^+}(q, {\bf t})\]
as in~(\ref{eq_factor}) above. 

\vspace{0.2in}

\noindent {\bf Acknowledgements} \ I would like to thank Aravind Asok, 
Kai Behrend, Tom Bridgeland, Jim Bryan, Alastair Craw, Victor Ginzburg, 
Dominic Joyce, Richard Kenyon, Rapha\"el Rouquier, 
Ed Segal, Yukinobu Toda and Ben Young for helpful
comments and correspondence, and the referee for his careful reading of the
manuscript. Special thanks are due to Alastair King for 
introducing me to quiver algebras and answering several questions; the main idea 
of the paper was formulated after a conversation with him. My research was 
partially supported by OTKA grant K61116. The final version of the paper was
prepared while I enjoyed the hospitality of the R\'enyi Institute of Mathematics 
in Budapest, Hungary, supported by the EU BudAlgGeo project.

\end{document}